\documentclass[11pt,reqno]{amsart}
\usepackage[T1]{fontenc}
\usepackage{lmodern,mathtools,amssymb,booktabs,array,microtype,needspace}
\usepackage[margin=1in]{geometry}
\usepackage[colorlinks=true,linkcolor=blue!50!black,citecolor=blue!50!black,urlcolor=blue!50!black]{hyperref}
\usepackage{xcolor}
\usepackage{graphicx}
\numberwithin{equation}{section}
\newtheorem{theorem}{Theorem}[section]
\newtheorem{proposition}[theorem]{Proposition}
\newtheorem{lemma}[theorem]{Lemma}
\newtheorem{corollary}[theorem]{Corollary}
\newtheorem{conjecture}[theorem]{Conjecture}
\theoremstyle{definition}
\newtheorem{remark}[theorem]{Remark}
\newcommand{\BR}{\operatorname{BR}}
\newcommand{\HH}{\mathcal H}
\newcommand{\cl}{\operatorname{cl}}

\newcommand{\cone}{\operatorname{cone}}
\newcommand{\FF}{\mathbb F}
\newcommand{\QQ}{\mathbb Q}
\newcommand{\T}{\mathsf T}

\title{Generic tensor products are not maximal}
\author{Mykhaylo Tyomkyn}
\address{Department of Applied Mathematics, Faculty of Mathematics and Physics, Charles University}
\email{tyomkyn@kam.mff.cuni.cz}
\thanks{Supported by GA\v{C}R grant 25-16847S}
\date{}
\begin{document}
\begin{abstract}
We construct tensor products of uniform matroids which are strictly freer than their generic linear tensor products. In particular, the generic tensor product of two copies of $U_{7,9}$ is not maximal. More generally, for every $a,b\ge2$, this holds for $U_{7,a+7}$ and $U_{7,b+7}$. The construction disproves the uniform tensor-product formulation of Mason's conjecture. Its dual is a relaxation of an explicit circuit-hyperplane of a generic birigidity matroid. By the author's representability theorem, the resulting matroids are non-representable over every field. We also obtain full-rank counterexamples to the Jackson--Tanigawa hyperconnectivity conjecture in every dimension $d\ge2$.
\end{abstract}
\maketitle

\section{Introduction}
We compare matroids by the weak order. That is, for matroids $M,N$ on the same finite ground set, write $M\preceq N$ if every independent set of $M$ is independent in $N$. We then say that $N$ is \emph{freer} than $M$ (or $M$ is \emph{less free} than $N$), and write $M\prec N$ if the two matroids are distinct. Consequently, we call a matroid \emph{maximum} or \emph{freest} in a class if it is freer than any other member of the class. We call a matroid  \emph{maximal} in a class if no other matroid in the class is freer.

Mason’s question on products of matroids dates from 1981~\cite{Mason}. Motivated by tensor products and powers of vector spaces, he asked whether the corresponding matroid constructions admit a freest member. For loopless matroids \(M_1,M_2\) on ground sets \(E_1,E_2\), a product is a matroid on \(E_1\times E_2\) whose restrictions to each row and column are naturally isomorphic to the corresponding factors, and in which \(F_1\times F_2\) is a flat whenever each \(F_i\) is a flat of \(M_i\). A product $T$ is a \emph{tensor product} if
\[
r_T(X_1\times X_2)=r_{M_1}(X_1)r_{M_2}(X_2)
\] 
for all $X_1\subseteq E_1$ and $X_2\subseteq E_2$.
Las Vergnas~\cite{LV} showed that a freest product need not exist in general, and that some pairs of matroids admit no tensor product at all.

For two uniform matroids, generic linear representations always provide a tensor product. The conjecture considered here asserts that this generic tensor product is their freest tensor product. This formulation was stated by Cruickshank, Jackson, Jordán and Tanigawa~\cite{CJJT} and, in its equivalent birigidity form, by Brakensiek, Chen, Dhar and Zhang~\cite{BCDZ}. In the latter paper the authors proved that the conjecture would imply a deterministic polynomial-time algorithm for testing independence in generic birigidity matroids, with corresponding applications to maximally recoverable tensor codes.

Write $U_{r,m}$ for the uniform matroid of rank $r$ on $[m]$.
Let $\FF$ be the rational function field over $\QQ$ in $mr+ns$
algebraically independent indeterminates. Arrange these indeterminates
as the coordinates of vectors $u_1,\ldots,u_m\in\FF^r$ and
$v_1,\ldots,v_n\in\FF^s$, using each indeterminate exactly once. The \emph{generic tensor product} $U_{r,m}\otimes U_{s,n}$
is the matroid on $[m]\times[n]$ in which a set
$S\subseteq[m]\times[n]$ is independent if and only if
the tensors $u_i\otimes v_j$, $(i,j)\in S$, are linearly
independent. This is a tensor product of $U_{r,m}$ and $U_{s,n}$ in the
matroidal sense defined above; see e.g.~\cite{BGILPS}.

\begin{conjecture}[Mason's conjecture, uniform tensor-product formulation]\label{conj:mason}
The generic tensor product of two uniform matroids is their freest tensor product.
\end{conjecture}

We disprove this conjecture by showing that the proposed maximum need not even be maximal.

\begin{theorem}\label{thm:main}
Let $k\ge2$, put $s=4k-1$, and let $a,b\ge2$. There is a tensor product $T_{a,b,k}$ of $U_{s,s+a}$ and $U_{s,s+b}$ such that
\begin{equation}\label{eq:main}
U_{s,s+a}\otimes U_{s,s+b}\prec T_{a,b,k}.
\end{equation}
The matroid $T_{a,b,k}$ is obtained by relaxing one circuit-hyperplane of the generic tensor product.
\end{theorem}

Taking $k=2$ gives the assertion for $U_{7,a+7}$ and $U_{7,b+7}$, and taking additionally $a=b=2$ gives two copies of $U_{7,9}$. The relaxed circuit-hyperplane in this last example has $49$ elements.

The proof is more transparent in the dual language of birigidity.
We construct a $32$-edge circuit-hyperplane of the generic birigidity
matroid $\BR_{2,2}(9,9)$, together with larger examples of the same
form, and use bipartite coning to increase the parameters $a$ and $b$.
Relaxing the resulting circuit-hyperplanes preserves every other
circuit and the total rank. This disproves the birigidity formulation
of Mason's conjecture~\cite{BCDZ, CJJT}, hence also the stronger
conjecture of Jackson and Tanigawa~\cite{JT},
and, by duality, Conjecture~\ref{conj:mason}.

The author proved~\cite{Tyomkyn} that the generic tensor product is the freest \emph{representable} tensor product of its two uniform factors. Equivalently, $\BR_{a,b}(m,n)$ is the freest representable abstract $(a,b)$-birigidity matroid on its ground set. Consequently, the counterexamples here are necessarily non-representable over every field.

Our second result concerns Kalai's generic hyperconnectivity
matroids. For integers $d\ge1$ and $n\ge d+1$, choose generic
vertex vectors $p_1,\ldots,p_n\in\FF^d$ over a rational
function field of characteristic $0$.
The generic $d$-hyperconnectivity matroid $\HH_d(K_n)$
is represented by the matrix whose column indexed by
$ij$, where $i<j$, has $p_j$ in the row block at $i$,
$-p_i$ in the row block at $j$, and zero elsewhere.
Its rank is
\[
r_d(n)=dn-\binom{d+1}{2},
\]
and every copy of $K_{d+2}$ and $K_{d+1,d+1}$ is a circuit.

For $d=2$, this matroid has a direct interpretation in
algebraic geometry: it is the algebraic matroid of the Grassmannian $\operatorname{Gr}(2,n)$
in its Pl\"ucker embedding.
Explicitly, writing $p_i=(x_i,y_i)^\T$, the Pl\"ucker
coordinates are
\[
q_{ij}=x_i y_j-x_j y_i
\qquad(1\le i<j\le n).
\]
An edge set $E$ is independent in $\HH_2(K_n)$ precisely
when the polynomials $q_{ij}$, $ij\in E$, are algebraically
independent over $\QQ$; see~\cite{Bernstein,TyomkynPlane}.

This interpretation also connects hyperconnectivity to
tropical geometry. Bernstein~\cite{Bernstein} used the
description of the tropical Grassmannian
$\operatorname{Gr}(2,n)$ in terms of phylogenetic trees
to characterize its independent edge sets:
a graph is $\HH_2$-independent if and only if it admits
an acyclic orientation with no alternating closed trail,
that is, no closed trail along which the edge directions
alternate.
The author~\cite{TyomkynPlane} gave a combinatorial proof
of the sufficiency direction and showed that the generic
$2$-hyperconnectivity matroid is independent of the field
characteristic.

Jackson and Tanigawa~\cite{JT} conjectured that
$\HH_d(K_n)$ is the unique maximal matroid on $E(K_n)$
in which every copy of $K_{d+2}$ and $K_{d+1,d+1}$
is a circuit.
Bulavka and Tancer~\cite{BT} disproved this conjecture
for $d\ge3$.
Our examples also cover $d=2$ and, in every dimension,
have full rigidity rank and are strictly freer than
$\HH_d(K_n)$.

\begin{theorem}\label{thm:hyper}
	Let $d,k\ge2$ be integers, and put $n=8k+d-1$.
	There is a matroid $M_{d,k}$ on $E(K_n)$ of rank $r_d(n)$
	such that
	\[
	\HH_d(K_n)\prec M_{d,k},
	\]
	and every copy of $K_{d+2}$ and $K_{d+1,d+1}$ is a circuit.
	The matroid $M_{d,k}$ is obtained by relaxing one
	circuit-hyperplane of $\HH_d(K_n)$ and is an abstract
	$d$-rigidity matroid.
	The complete bipartite condition is non-vacuous whenever
	$8k\ge d+3$.
\end{theorem}

Taking $d=k=2$ gives a counterexample on $17$ vertices.
Thus the algebraic matroid of $\operatorname{Gr}(2,17)$ in its Pl\"ucker embedding is not maximal among matroids
in which every $K_4$ and $K_{3,3}$ is a circuit.
The construction and proof are given in
Section~\ref{sec:hypercon}.

\subsection*{Statement on AI use}
The counterexamples were found by ChatGPT 6 in the course of the following interaction. Jackson and Tanigawa~\cite{JT} asked if $\BR_{2,2}(m,n)$ matroids have the $K_{3,3}$-covering property. That is, if every cyclic flat is a union of copies of $K_{3,3}$. If true, then, as was shown in~\cite{JT}, this would imply the Jackson-Tanigawa conjecture, and a fortiori the Mason conjecture for all $\BR_{2,2}(m,n)$. I was trying to prove a weaker statement, namely, that any $\BR_{2,2}(m,n)$-circuit that is also a flat must be a copy of $K_{3,3}$. I prompted ChatGPT multiple ideas, and at first it made good partial progress. For instance, it handled the cases when all degree $3$ vertices are on the same side, and when two degree $3$ vertices are adjacent. This apparently narrowed down its search space sufficiently that it found the $9\times 9$-counterexample. The rest was routine.

I have revised and rewritten the proofs, and take responsibility for the content and wording of the paper.

\section{Preliminaries and the construction}
Let $L,R$ be disjoint sets of sizes $m,n$. We will refer to members of $L$ and $R$ as the `left' and `right' vertices, respectively. Identify $L\times R$ with the edges of $K_{m,n}$. Fix positive integers $a,b$ with $m\ge a+1$ and $n\ge b+1$. An \emph{abstract $(a,b)$-birigidity matroid} is a matroid on $L\times R$ of rank
\begin{equation}\label{eq:ambient-rank}
r_{a,b}(m,n)=bm+an-ab
\end{equation}
in which every $A\times B$, with $|A|=a+1$ and $|B|=b+1$, is a circuit. 

For a (column) matrix representation $A$ of a matroid, we call vectors in $\ker A$ \emph{stresses} and vectors in $\ker A^\T$ \emph{motions}.
For a set $F$ of columns, its stress space is the kernel
of the submatrix consisting of those columns.
Thus $F$ is a circuit precisely when its stress space is
one-dimensional and a nonzero stress has no zero coordinates.

One central example of an abstract birigidity matroid is the \emph{generic birigidity matroid} $\BR_{a,b}(m,n)$, which has the following representation. Choose generic vectors $u_i\in\FF^a$ for $i\in L$ and $v_j\in\FF^b$ for $j\in R$, over a rational function field $\FF$ of characteristic $0$. The column indexed by $ij$ has $v_j$ in the $b$ rows belonging to $i$, $-u_i$ in the $a$ rows belonging to $j$, and zero elsewhere. We use the signed convention so that, when $a=b=d$, the
birigidity representation agrees with the restriction of (the column form of) the $d$-hyperconnectivity representation to a bipartite graph. Changing all signs in the right vertex blocks gives the usual unsigned convention.

Write $U$ and $V$ for the matrices whose rows are $u_i^{\T}$
and $v_j^{\T}$. A $\BR_{a,b}(m,n)$-stress of a bipartite
graph $G\subseteq L\times R$ is equivalently a matrix $W$,
supported on $E(G)$, satisfying
\begin{equation}\label{eq:stress}
	WV=0,\qquad W^{\T}U=0.
\end{equation}
Every copy of $K_{a+1,b+1}$ is a circuit in
$\BR_{a,b}(m,n)$~\cite{KNN}.

A $\BR_{a,b}(m,n)$-motion of $G$ consists of vectors
$\dot u_i\in\FF^b$ and $\dot v_j\in\FF^a$ satisfying
\begin{equation*}
	v_j^{\T}\dot u_i-u_i^{\T}\dot v_j=0
	\qquad(ij\in E(G)).
\end{equation*}
For the complete bipartite graph, these motions are exactly
$\dot u_i=A u_i$, $\dot v_j=A^{\T}v_j$, where $A$ ranges
over all $b\times a$ matrices~\cite{KNN}.
Their space has dimension $ab$, giving~\eqref{eq:ambient-rank}
by rank--nullity.

The duality between generic tensor products and generic birigidity was established in~\cite{BDGGL}. Its abstract version is described in~\cite{CJJT}: a matroid $N$ on $[m]\times[n]$ is an abstract $(a,b)$-birigidity matroid if and only if $N^*$ is a tensor product of $U_{m-a,m}$ and $U_{n-b,n}$. In particular,
\begin{equation}\label{eq:duality}
\BR_{a,b}(m,n)^*=U_{m-a,m}\otimes U_{n-b,n}.
\end{equation}

A \emph{circuit-hyperplane} of a rank-$r$ matroid $M$ is a set $F$ which is both a circuit and a hyperplane. In particular, $|F|=r$, $r_M(F)=r-1$, and $\cl_M(F)=F$. We recall the circuit-hyperplane relaxation method, which is well-known.

\begin{lemma}\label{lem:relaxation}
If $F$ is a circuit-hyperplane of $M$, adjoining $F$ to the bases of $M$ defines a matroid $M^F$, called its relaxation. Its rank function differs from that of $M$ only at $F$, where it increases by $1$. Every circuit of $M$ other than $F$ remains a circuit of $M^F$. Moreover, $(M^F)^*$ is the relaxation of the circuit-hyperplane $E(M)\setminus F$ in $M^*$.
\end{lemma}
\begin{proof}
We check basis exchange for an old basis $B$ and the new basis $F$. For $x\in F\setminus B$, the set $F-x$ spans in $M$ the hyperplane $F$. The basis $B$ contains an element $y$ outside the $M$-hyperplane $F$, and therefore $F-x+y$ is an old basis for an element $y\in B\setminus F$. Conversely, fix $x\in B\setminus F$. If $F$ is not contained in $\cl_M(B-x)$, choose $y\in F\setminus\cl_M(B-x)$; then $B-x+y$ is an old basis. Otherwise, $\cl_M(B-x)=F$, and hence $B-x\subset F$. The unique element $y\in F\setminus B$ gives $B-x+y=F$.

Every proper subset of $F$ was already independent, and every proper superset of $F$ already had rank $r$. Thus only the rank of $F$ changes. A circuit different from $F$ cannot contain $F$, so it remains a circuit. The last assertion follows from the duality between hyperplanes and circuits.
\end{proof}
Fix an integer $k\ge2$. In the construction below, $i$ ranges over $\mathbb Z/k\mathbb Z$. Put
\begin{align*}
L_k=\{a_*\}\cup\{a_i^1,a_i^2,x_i^1,x_i^2:i\in\mathbb Z/k\mathbb Z\} \ \ \text{ and } \ \ R_k=\{b_*\}\cup\{b_i^1,b_i^2,y_i^1,y_i^2:i\in\mathbb Z/k\mathbb Z\}.
\end{align*}
Define $G_k$ by the following edges, for each $i$ and $\sigma\in\{1,2\}$:
\begin{equation*}
\begin{gathered}
x_i^\sigma b_*,\quad x_i^\sigma b_i^1,\quad x_i^\sigma b_i^2,\quad a_*y_i^\sigma,\quad a_i^1y_i^\sigma,\quad a_i^2y_i^\sigma, \quad a_i^\sigma b_i^\sigma,\quad a_i^\sigma b_{i-1}^\sigma.
\end{gathered}
\end{equation*}
The first six edge types above give $2k$ copies of $K_{2,3}$, which we call the \emph{lateral blocks}. The edges of the last two types are the \emph{connectors}. There are $4k+1$ vertices in each class and $16k$ edges. For $k=2$, the construction consists of four lateral blocks and two $K_{2,2}$ connector blocks.

\begin{figure}[!b]
\centering
\begin{minipage}[c]{0.48\linewidth}
\includegraphics[width=\linewidth]{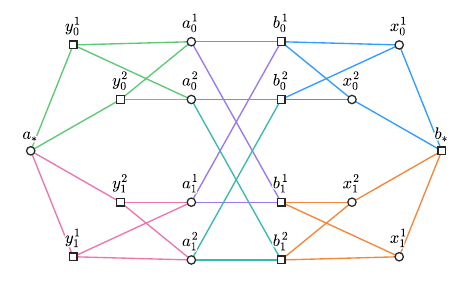}
\end{minipage}\hfill
\begin{minipage}[c]{0.48\linewidth}
\includegraphics[width=\linewidth]{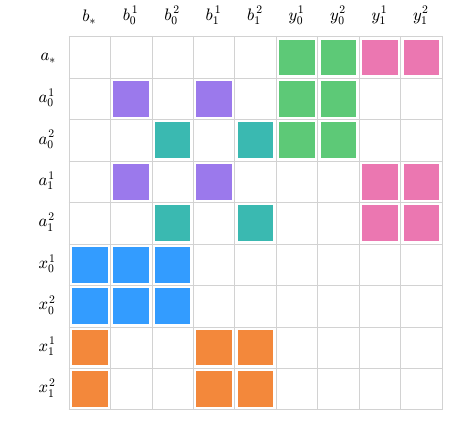}
\end{minipage}
\caption{The graph $G_2$ and its adjacency matrix. Circles represent $L_2$
and squares represent $R_2$; matrix rows and columns follow these classes.
Each lateral block has its own colour, and the two connector $4$-cycles
are purple and turquoise. Filled cells represent edges of the same colour.}
\label{fig:c2}
\end{figure}

\begin{figure}[!t]
\centering
\begin{minipage}[c]{0.48\linewidth}
\includegraphics[width=\linewidth]{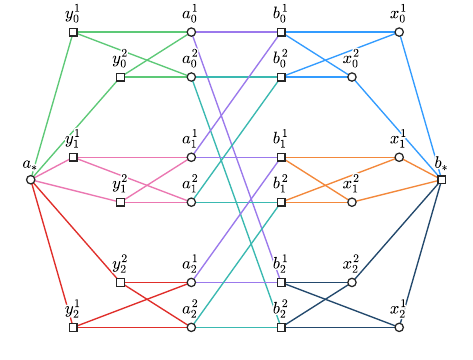}
\end{minipage}\hfill
\begin{minipage}[c]{0.48\linewidth}
\includegraphics[width=\linewidth]{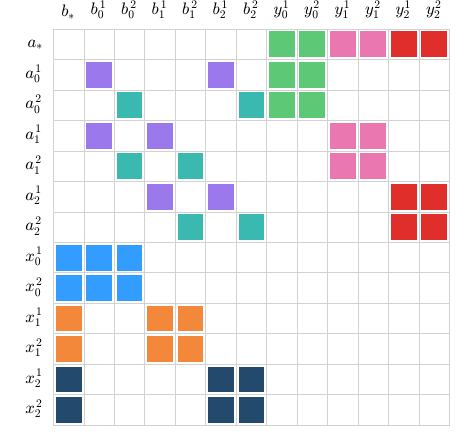}
\end{minipage}
\caption{The graph $G_3$ and its adjacency matrix, with the same conventions
as in Figure~\ref{fig:c2}. The six lateral blocks have distinct colours;
the purple and turquoise connectors each form a $6$-cycle.}
\label{fig:c3}
\end{figure}

The graphs $G_2$ and $G_3$, together with their adjacency matrices,
are shown in Figures~\ref{fig:c2} and~\ref{fig:c3}.

\section{The proof}
\begin{proposition}\label{prop:seed}
	For every $k\ge2$, the edge set of $G_k$ is a circuit-hyperplane
	of $\BR_{2,2}(4k+1,4k+1)$.
\end{proposition}
We first bound the dimension of the stress space, then use
an explicit motion to establish the rank, circuit and flatness
assertions. Throughout this section, $i,j$ range over $\mathbb Z/k\mathbb Z$
and $\sigma,\tau$ over $\{1,2\}$, unless stated otherwise.

\subsection{The stress space}
Write $p_v\in\FF^2$ for the generic vector assigned to vertex $v$,
using the same notation for both vertex classes. The represented
matroid is unchanged by a common invertible linear transformation
of these vectors or by multiplying individual vectors by nonzero
scalars; see~\cite[Lemma~2.1]{CRS}. We use these operations to
normalize the vectors. Since $p_{a_*},p_{b_*}$ are
independent, a common change of coordinates sends them to $(1,0)^\T$ and $(0,1)^\T$ respectively.
Individual rescalings then give
\begin{equation*}
	\begin{gathered}
		p_{a_*}=(1,0)^\T,\qquad p_{b_*}=(0,1)^\T,\\
		p_{a_i^1}=(w_i,1)^\T,\quad p_{a_i^2}=(z_i,1)^\T,\qquad
		p_{b_i^1}=(1,u_i)^\T,\quad p_{b_i^2}=(1,v_i)^\T.
	\end{gathered}
\end{equation*}
The parameters $w_i,z_i,u_i,v_i$ and the coordinates at the vertices
$x_i^\sigma,y_i^\sigma$ remain algebraically independent, as
the original coordinates can be recovered rationally from the
normalized coordinates, the original basis vectors
$p_{a_*},p_{b_*}$, and the scaling factors.

Put
\[
d_i=w_i-z_i,\qquad e_i=u_i-v_i.
\]
These quantities are nonzero by genericity.

Fix an arbitrary stress $\omega$ on $G_k$, and write
$\omega(ab)$ for the coefficient of edge $ab$. By~\eqref{eq:stress},
the stress equation at any vertex $v$ is
\[
\sum_{t\in N(v)}\omega(vt)p_t=0,
\]
where $\omega(vt)=\omega(tv)$ denotes the coefficient of the same
edge. 

We first determine the connector coefficients. Write
$A_i^\sigma=\omega(a_i^\sigma b_i^\sigma)$ and
$B_i^\sigma=\omega(a_i^\sigma b_{i-1}^\sigma)$. Since the first coordinates of $p_{b_i^1},p_{b_i^2},p_{b_*}$
are $1,1,0$, respectively, the stress equation at $x_i^\sigma$
forces the coefficients of $x_i^\sigma b_i^1$ and
$x_i^\sigma b_i^2$ to be opposites. Similarly, the second coordinate at $y_i^\sigma$ forces the coefficients of its edges to $a_i^1,a_i^2$ to be
opposites. Adding the stress equations at $a_i^1,a_i^2$, or at
$b_i^1,b_i^2$, therefore cancels the lateral contributions and gives
\begin{align}\label{eq:coordinateseqn}
	A_i^1p_{b_i^1}+A_i^2p_{b_i^2}
	+B_i^1p_{b_{i-1}^1}+B_i^2p_{b_{i-1}^2}&=0, \ \text{ and } \nonumber \\
	A_i^1p_{a_i^1}+A_i^2p_{a_i^2}
	+B_{i+1}^1p_{a_{i+1}^1}+B_{i+1}^2p_{a_{i+1}^2}&=0.
\end{align}
Taking the first coordinate of the first equation and the second
coordinate of the second shows that
\[
A_i^1+A_i^2=-(B_i^1+B_i^2)=-(B_{i+1}^1+B_{i+1}^2).
\]
Hence there is a common scalar $\rho$ such that
$A_i^1+A_i^2=\rho$ and $B_i^1+B_i^2=-\rho$ for every $i$.
Define
\[
q_i=d_i e_i A_i^1,\qquad r_i=d_i e_{i-1}B_i^1.
\]
The remaining coordinate equations of~\eqref{eq:coordinateseqn}, multiplied by $d_i$ and $e_i$,
respectively, become
\[
\begin{aligned}
	q_i+r_i&=-\rho d_i(v_i-v_{i-1}),\\
	q_i+r_{i+1}&=-\rho e_i(z_i-z_{i+1}).
\end{aligned}
\]
Put
\[
P=\sum_i w_i(v_i-v_{i-1}) \ \text{ and } \ 
Q=\sum_i z_i(u_i-u_{i-1}).
\]
Subtracting the second equation from the first gives
\[
r_i-r_{i+1}
=\rho\bigl(e_i(z_i-z_{i+1})-d_i(v_i-v_{i-1})\bigr).
\]
When we sum over all cyclic indices, the left-hand side
vanishes. To evaluate the right-hand side, reindex the sum
containing $z_{i+1}$ to obtain
\[
\sum_i e_i(z_i-z_{i+1})
=\sum_i z_i(e_i-e_{i-1}).
\]
Using $e_i=u_i-v_i$ and $d_i=w_i-z_i$, we therefore have
\[
\begin{aligned}
	&\sum_i\bigl(e_i(z_i-z_{i+1})-d_i(v_i-v_{i-1})\bigr)\\
	&\quad=
	\sum_i z_i(u_i-u_{i-1})
	-\sum_i z_i(v_i-v_{i-1})
	-\sum_i(w_i-z_i)(v_i-v_{i-1})\\
	&\quad=
	\sum_i z_i(u_i-u_{i-1})
	-\sum_i w_i(v_i-v_{i-1})
	=Q-P.
\end{aligned}
\]
Thus summing the equations for $r_i-r_{i+1}$ gives
$\rho(Q-P)=0$.
For $k\ge2$, the monomial $z_0u_0$ has coefficient $1$ in $Q$
and does not occur in $P$. Thus $Q-P\ne0$. It follows that
$\rho=0$, so $q_i=-r_i=-r_{i+1}$ for every $i$.

The connector coefficients therefore depend on at most one scalar.
We claim that they uniquely determine all lateral coefficients.
Indeed, the stress equation at $b_i^1$ is
\[
\omega(x_i^1b_i^1)p_{x_i^1}
+\omega(x_i^2b_i^1)p_{x_i^2}
=-A_i^1p_{a_i^1}-B_{i+1}^1p_{a_{i+1}^1}.
\]
Its right-hand side is determined by the connector coefficients.
Since $p_{x_i^1},p_{x_i^2}$ form a basis of $\FF^2$, this equation
uniquely determines $\omega(x_i^1b_i^1)$ and
$\omega(x_i^2b_i^1)$. The two coordinate equations at each
$x_i^\sigma$ then give
\[
\begin{aligned}
	\omega(x_i^\sigma b_i^2)
	&=-\omega(x_i^\sigma b_i^1) \ \ \text{ and }\\
	\omega(x_i^\sigma b_*)
	&=-e_i\,\omega(x_i^\sigma b_i^1).
\end{aligned}
\]
Thus all coefficients in these lateral blocks are determined.

Similarly, the stress equation at $a_i^1$, together with the
linear independence of $p_{y_i^1},p_{y_i^2}$, uniquely determines
$\omega(a_i^1y_i^1)$ and $\omega(a_i^1y_i^2)$. The equations at
the $y_i^\sigma$ then determine the coefficients of their edges
to $a_i^2$ and $a_*$. Consequently, the entire stress space has
dimension at most $1$.

\subsection{An explicit motion}
To show that the above stress space is non-trivial, we exhibit a motion. To do so, we first choose the motion vectors at the endpoints of the connectors,
then extend it to the lateral blocks.

Let $\lambda,\mu,\beta_i,\gamma_i\in\FF$ be scalars to be
determined. Recall that $d_i=w_i-z_i$. Set
\[
\dot p_{b_i^\sigma}
=(\lambda u_i-\mu v_i,\beta_i)^\T,
\]
\[
\dot p_{a_i^1}=(\gamma_i,\lambda d_i)^\T
\qquad\text{and}\qquad
\dot p_{a_i^2}=(\gamma_i,\mu d_i)^\T.
\]
For any pair $ab\in L_k\times R_k$, write
\[
\varepsilon(ab)=p_b^\T\dot p_a-p_a^\T\dot p_b.
\]
Using $d_i=w_i-z_i$, we obtain
\[
\begin{aligned}
	\varepsilon(a_i^1b_j^1)
	&=\gamma_i-\beta_j+\lambda d_i u_j
	-w_i(\lambda u_j-\mu v_j)\\
	&=\gamma_i-\beta_j-\lambda z_i u_j+\mu w_i v_j,
\end{aligned}
\]
and
\[
\begin{aligned}
	\varepsilon(a_i^2b_j^2)
	&=\gamma_i-\beta_j+\mu d_i v_j
	-z_i(\lambda u_j-\mu v_j)\\
	&=\gamma_i-\beta_j-\lambda z_i u_j+\mu w_i v_j.
\end{aligned}
\]
Thus both types of connector impose the same condition:
\[
\gamma_i-\beta_j
=\lambda z_i u_j-\mu w_i v_j
\qquad\text{for }j\in\{i,i-1\}.
\]
Subtracting the two conditions gives
\[
\beta_{i-1}-\beta_i
=\lambda z_i(u_i-u_{i-1})
-\mu w_i(v_i-v_{i-1}).
\]
This recurrence is consistent around the cycle precisely when
the sum of its right-hand sides vanishes, that is, when
\[
\lambda Q-\mu P=0.
\]
We therefore take $\lambda=P$ and $\mu=Q$.

Set $\beta_0=0$ and define $\beta_1,\ldots,\beta_{k-1}$
successively by
\begin{equation}\label{eq:recurrence}
	\beta_{i-1}-\beta_i
	=Pz_i(u_i-u_{i-1})-Qw_i(v_i-v_{i-1}).
\end{equation}
The sum of the right-hand sides is $PQ-QP=0$, so the equation
with $i=0$ holds as well. Define
\[
\gamma_i=\beta_i+Pz_i u_i-Qw_i v_i.
\]
Every connector equation now holds. More generally, the
calculation above gives
\[
\varepsilon(a_i^1b_j^1)
=\varepsilon(a_i^2b_j^2)
=\beta_i-\beta_j
+Pz_i(u_i-u_j)-Qw_i(v_i-v_j).
\]

To extend this assignment to the lateral blocks, put
\[
h=P-Q,\qquad
\alpha_i=Pu_i-Qv_i,\qquad\text{and}\qquad
\delta_i=Qw_i-Pz_i,
\]
and define
\[
M_i=
\begin{pmatrix}
	\alpha_i&\beta_i\\
	0&0
\end{pmatrix}
\qquad\text{and}\qquad
N_i=
\begin{pmatrix}
	0&\gamma_i\\
	h&\delta_i
\end{pmatrix}.
\]
These matrices give the already prescribed motions at the
connector endpoints:
\[
M_i^\T p_{b_i^\sigma}
=(\alpha_i,\beta_i)^\T
=\dot p_{b_i^\sigma},
\]
\[
N_i p_{a_i^1}
=(\gamma_i,Pd_i)^\T
=\dot p_{a_i^1}
\qquad\text{and}\qquad
N_i p_{a_i^2}
=(\gamma_i,Qd_i)^\T
=\dot p_{a_i^2}.
\]
Here the last two identities use
\[
hw_i+\delta_i=Pd_i
\qquad\text{and}\qquad
hz_i+\delta_i=Qd_i.
\]
At the remaining vertices, set
\[
\dot p_{x_i^\sigma}=M_i p_{x_i^\sigma},
\qquad
\dot p_{y_i^\sigma}=N_i^\T p_{y_i^\sigma},
\]
\[
\dot p_{b_*}=(0,0)^\T
\qquad\text{and}\qquad
\dot p_{a_*}=(0,h)^\T,
\]
and note that
\[
M_i^\T p_{b_*}=(0,0)^\T
\qquad\text{and}\qquad
N_i p_{a_*}=(0,h)^\T
\]
for every $i$.

For any $2\times2$ matrix $A$, assigning
$\dot p_a=Ap_a$ and $\dot p_b=A^\T p_b$ gives
\[
p_b^\T\dot p_a-p_a^\T\dot p_b
=p_b^\T A p_a-p_a^\T A^\T p_b=0.
\]
Applying this observation with $A=M_i$ to every edge incident
to $x_i^1$ or $x_i^2$, and with $A=N_i$ to every edge incident
to $y_i^1$ or $y_i^2$, shows that all lateral edge equations
are satisfied. Thus the assignment is a motion of $G_k$.

Our motion has
\[
\varepsilon(a_*b_*)=h=P-Q\ne0,
\]
and so the column of $a_*b_*$ is outside the span of the columns
of $G_k$. Since the complete bipartite host has rank $16k$,
we obtain $r_{\BR_{2,2}}(G_k)\le16k-1$. There are $16k$ edges,
so a nonzero stress exists. Together with the upper bound
proved above, this shows that the stress space is exactly
one-dimensional, and hence
\begin{equation}\label{eq:seed-rank}
	r_{\BR_{2,2}}(G_k)=16k-1.
\end{equation}

We now prove that a nonzero stress has full support.
Its common value of $r_i$ must be nonzero: otherwise all
connector coefficients, and hence all lateral coefficients,
would vanish. The formulas above therefore show that every
connector coefficient is nonzero.

The support of this stress is the unique circuit $C$ contained
in $G_k$. Generic birigidity is invariant under relabelling
within either vertex class, so every graph automorphism of $G_k$ preserving the two vertex classes preserves $C$. In particular, this holds for the
interchanges $x_i^1\leftrightarrow x_i^2$ and $y_i^1\leftrightarrow y_i^2$.

Every vertex of a circuit has degree at least $3$, since
any two generic neighbour vectors are linearly independent.
At $b_i^1$, both connector edges belong to $C$, so at least
one of its edges to $x_i^1,x_i^2$ must also belong to $C$.
Interchanging these two vertices shows that both edges belong
to $C$. Each $x_i^\sigma$ now belongs to $C$ and has degree
$3$ in $G_k$, so all three of its incident edges belong
to $C$. The same argument at $a_i^1$, using the interchange
of $y_i^1,y_i^2$, includes every edge of the other lateral
blocks. Thus $C=G_k$, proving the circuit assertion.

\subsection{Flatness}
The value $\varepsilon(ab)$ is the evaluation of the edge column
against our motion. It vanishes on every existing edge and hence
on the span of their columns. Consequently, to prove flatness it
suffices to show that $\varepsilon(ab)\ne0$ for every missing
crossedge. Since $\beta_0=0$, the recurrence~\eqref{eq:recurrence}
expresses every $\beta_i$ as a polynomial in the parameters
$w_\ell,z_\ell,u_\ell,v_\ell$,
$\ell\in\mathbb Z/k\mathbb Z$.
The same is true of $h,\alpha_i,\gamma_i,\delta_i$.
Consequently, every missing-edge evaluation is a polynomial
in these parameters and the coordinates at the vertices
$x_i^\sigma,y_i^\sigma$.
Since these quantities are algebraically independent,
it suffices to show that each such polynomial is not
identically zero.

For each $i,j$, put
\[
T_{ij}=\beta_i-\beta_j
+Pz_i(u_i-u_j)-Qw_i(v_i-v_j).
\]
Thus $\varepsilon(a_i^1b_j^1)=\varepsilon(a_i^2b_j^2)=T_{ij}$.
The pairs with $j=i$ or $j=i-1$ are existing edges. For $k=2$, every pair $a_i^\sigma b_j^\sigma$ is an existing edge. For $k\geq 3$, for a missing pair $a_i^\sigma b_j^\sigma$, we may choose the indices so that $j+2\le i\le j+k-1$. Summing
\eqref{eq:recurrence} along this interval gives
\[
\beta_i-\beta_j
=-P\sum_{q=j+1}^{i}z_q(u_q-u_{q-1})
+Q\sum_{q=j+1}^{i}w_q(v_q-v_{q-1}).
\]
Define
\[
A_{ij}(z,u)
=\sum_{q=j+1}^{i-1}(z_i-z_q)(u_q-u_{q-1}),
\]
and define $A_{ij}(w,v)$ analogously. Using
$u_i-u_j=\sum_{q=j+1}^{i}(u_q-u_{q-1})$ and the corresponding
identity for $v$, we obtain
\[
T_{ij}=PA_{ij}(z,u)-QA_{ij}(w,v).
\]
The polynomial $A_{ij}(w,v)$ is nonzero: its coefficient of
$w_{j+1}$ is $-(v_{j+1}-v_j)$. Also, $A_{ij}(z,u)$ does not involve
$z_j$, and $P$ is independent of all the $z$-variables.
The coefficient of $z_j$ in $T_{ij}$ is therefore
\[
-(u_j-u_{j-1})A_{ij}(w,v),
\]
which is nonzero. Thus every missing pair $a_i^\sigma b_j^\sigma$ has
nonzero evaluation.

For the pairs $a_i^1b_j^2$ and $a_i^2b_j^1$,
direct substitution gives
\begin{align*}
	\varepsilon(a_i^1b_j^2)
	&=T_{ij}+P(w_i-z_i)(v_j-u_j)\ \ \text{ and } \\
	\varepsilon(a_i^2b_j^1)
	&=T_{ij}+Q(w_i-z_i)(u_j-v_j).
\end{align*}
To see that the first polynomial is nonzero, set all $z_q$ and
$u_q$ equal to zero. Then $Q=0$, and~\eqref{eq:recurrence} with
$\beta_0=0$ gives $\beta_q=0$ for every $q$. Hence $T_{ij}=0$,
and the first expression becomes $Pw_i v_j$, a nonzero polynomial
in the remaining variables. Similarly, setting all $w_q$ and
$v_q$ equal to zero reduces the second expression to $-Qz_i u_j$,
which is also nonzero. These substitutions are used only to test
whether the polynomial expressions vanish identically.
Thus every pair of the form $a_i^1b_j^2$ or $a_i^2b_j^1$
has nonzero evaluation.

It remains to consider crossedges involving $a_*,b_*$ or vertices
$x_i^\sigma,y_j^\tau$. For $x=x_i^\sigma$ and $y=y_j^\tau$, write
$p_x=(s,t)^\T$ and $p_y=(r,\ell)^\T$. Their motions are
\[
\dot p_x=(\alpha_i s+\beta_i t,0)^\T,\qquad
\dot p_y=(h\ell,\gamma_j r+\delta_j\ell)^\T.
\]
So, for example,
\[
\varepsilon(x_i^\sigma y_j^\tau)
=\alpha_i sr+(\beta_i-\gamma_j)tr-hs\ell-\delta_jt\ell.
\]
The remaining missing pairs give the following evaluations or
coefficients. Coefficients are taken with respect to the coordinates
$s,t,r,\ell$, which are algebraically independent of the parameters
$w_i,z_i,u_i,v_i$.
\begin{center}
	\renewcommand{\arraystretch}{1.22}
	\begin{tabular}{@{}ll@{}}
		\toprule
		Missing pair&Evaluation, or indicated coefficient\\
		\midrule
		$a_*b_*$&$P-Q$\\
		$a_*b_i^1,\ a_*b_i^2$&$-Q(u_i-v_i),\ -P(u_i-v_i)$\\
		$a_i^1b_*,\ a_i^2b_*$&$P(w_i-z_i),\ Q(w_i-z_i)$\\
		$a_i^\sigma y_j^\tau$, $i\ne j$
		&coefficient of $\ell$: $\delta_i-\delta_j$\\
		$x_i^\sigma b_j^\tau$, $i\ne j$
		&coefficient of $s$: $\alpha_i-\alpha_j$\\
		$x_i^\sigma y_j^\tau$
		&coefficient of $s\ell$: $-(P-Q)$\\
		\bottomrule
	\end{tabular}
\end{center}
The monomials $w_0v_0$ in $P$ and $z_0u_0$ in $Q$ both
have coefficient $1$, so $P\ne0$ and $Q\ne0$.
We have already shown that $P-Q\ne0$, while
$u_i-v_i\ne0$ and $w_i-z_i\ne0$ by genericity. For $i\ne j$, we have
\begin{align*}
	\delta_i-\delta_j
	&=Q(w_i-w_j)-P(z_i-z_j),\\
	\alpha_i-\alpha_j
	&=P(u_i-u_j)-Q(v_i-v_j).
\end{align*}
The coefficient of $u_0$ in the first polynomial is
$(z_0-z_1)(w_i-w_j)$, while the coefficient of $z_0$ in the second
is $-(u_0-u_{k-1})(v_i-v_j)$. Both are nonzero.
Thus every expression used in the table is nonzero.

Since the coordinates at the vertices $x_i^\sigma,y_j^\tau$ are
algebraically independent of the other parameters, a nonzero
coefficient in these coordinates guarantees a nonzero evaluation.
The table, together with the pairs $a_i^\sigma b_j^\tau$ treated
above, covers every missing crossedge. Hence none of their columns
belongs to the span of the columns of $G_k$, proving flatness.

By~\eqref{eq:seed-rank}, this flat has rank $16k-1$, whereas
\eqref{eq:ambient-rank} gives ambient rank $16k$. It is therefore
a hyperplane. Together with the circuit assertion, this proves
Proposition~\ref{prop:seed}.

\section{Increasing the birigidity parameters}

The following is a circuit-hyperplane version of the
bipartite coning lemma of Kalai, Nevo and
Novik~\cite[Lemma~3.12]{KNN}; see also~\cite{BDGGL}.
We include a proof.
\begin{lemma}[Bipartite coning]\label{lem:bicone}
	Let $a,b,m,n$ be positive integers with $m\ge a+1$ and
	$n\ge b+1$. Suppose that $F\subseteq K_{m,n}$ uses every
	vertex and is a circuit-hyperplane of $\BR_{a,b}(m,n)$.
	Add one new left vertex adjacent to every right vertex,
	obtaining $\widehat F$. Then $\widehat F$ is a
	circuit-hyperplane of $\BR_{a+1,b}(m+1,n)$.
	The analogous assertion holds on the right, with $b$
	increased by $1$.
\end{lemma}
\begin{proof}
Consider the generic representation of $\BR_{a+1,b}(m+1,n)$,
with vectors in $\FF^{a+1}$ assigned to the left vertices
and vectors in $\FF^b$ assigned to the right vertices.
Label the old left and right vertices by $[m]$ and $[n]$,
respectively. An invertible change of coordinates applied
to all left vectors preserves the represented matroid,
so we may take the vector at the new left vertex to be
$(0,\ldots,0,1)^\T$.
Write the old left vectors as $(u_i^\T,t_i)^\T$, where
$u_i\in\FF^a$ and $i\in[m]$.
	Write $v_j\in\FF^b$ for the right vectors, $j\in[n]$.
	The coordinates of the $u_i,v_j$ and the scalars $t_i$
	are generic. Put
	\[
	t=(t_1,\ldots,t_m)^{\T}.
	\]
Let $U$ be the $m\times a$ matrix with $i$th row $u_i^\T$,
and let $V$ be the $n\times b$ matrix with $j$th row $v_j^\T$.
	
	For a graph $G$ on the old vertex classes, let $\widehat G$
	be obtained by adding the new left vertex and all its $n$
	incident edges. We record stress coefficients in a matrix whose rows and
	columns correspond to the left and right vertices, respectively,
	putting zero in positions corresponding to missing edges.
	For $\widehat G$, place the new left vertex last.
	Then the coefficient matrix of every stress of $\widehat G$
	must have the form
	\[
	\widehat W=
	\begin{pmatrix}
		W\\-t^{\T}W
	\end{pmatrix},
	\qquad
	WV=0,\quad W^{\T}U=0,
	\]
	where $W$ is supported on $E(G)$.
	Indeed, the last-coordinate equations at the right vertices
	force the new row to be $-t^{\T}W$. The remaining equations
	at the old vertices are exactly~\eqref{eq:stress}, and the
	equation at the new vertex follows from $WV=0$.
	Thus the old and new stress spaces have the same dimension.
	Since $n$ edges were added,
	\begin{equation}\label{eq:bicone-rank}
		r_{\BR_{a+1,b}}(\widehat G)
		=r_{\BR_{a,b}}(G)+n.
	\end{equation}
	
	For $G=F$, choose a nonzero old stress $W$. It has full
	support, and every column of $W$ is nonzero because every
	right vertex occurs in $F$. The scalars $t_i$ are
	algebraically independent of the old coordinates, so each
	new coefficient
	\[
	-\sum_{i=1}^m t_iW_{ij}
	\qquad(j\in[n])
	\]
	is nonzero. Hence the new stress space is one-dimensional
	and its nonzero stresses have full support, proving that
	$\widehat F$ is a circuit.
	
	Every missing edge of $\widehat F$ is an old missing edge.
	Applying~\eqref{eq:bicone-rank} to $F$ and to $F+f$ for
	each such edge $f$ shows that flatness is preserved.
	Finally,
	\[
	r_{a+1,b}(m+1,n)=r_{a,b}(m,n)+n,
	\]
	so $\widehat F$ has rank one less than its ambient matroid.
	It is therefore a hyperplane. Transposition gives the
	right-hand version.
\end{proof}

For integers $a,b,k\ge2$, let $G_k^{a,b}$ be obtained from
$G_k$ by adding $a-2$ universal left vertices and $b-2$
universal right vertices. In particular, every new left
vertex is joined to every new right vertex. Its vertex
class sizes and edge count are
\begin{align}
	m&=4k+a-1,\qquad n=4k+b-1,
	\label{eq:coned-sizes}\\
	|E(G_k^{a,b})|
	&=16k+(4k+1)(a+b-4)+(a-2)(b-2)\notag\\
	&=ab+(4k-1)(a+b)=r_{a,b}(m,n).\notag
\end{align}

\begin{theorem}\label{thm:birigidity}
	For every choice of integers $a,b,k\ge2$, the graph
	$G_k^{a,b}$ is a circuit-hyperplane of $\BR_{a,b}(m,n)$,
	with $m,n$ as in~\eqref{eq:coned-sizes}.
	Its relaxation $N_{a,b,k}$ is an abstract
	$(a,b)$-birigidity matroid and
	\[
	\BR_{a,b}(m,n)\prec N_{a,b,k}.
	\]
\end{theorem}

\begin{proof}
	Proposition~\ref{prop:seed} and repeated application of
	Lemma~\ref{lem:bicone} give the circuit-hyperplane assertion.
	The graph uses all vertices and
	\[
	m-(a+1)=n-(b+1)=4k-2>0,
	\]
	so it is not a $K_{a+1,b+1}$.
	By Lemma~\ref{lem:relaxation}, every prescribed
	$K_{a+1,b+1}$ remains a circuit after relaxation.
	The total rank is unchanged. These are exactly the defining
	conditions for abstract birigidity, and the relaxed set
	witnesses strictness.
\end{proof}

\begin{proof}[Proof of Theorem~\ref{thm:main}]
	Fix integers $a,b,k\ge2$, put $s=4k-1$, and take
	$T_{a,b,k}=N_{a,b,k}^*$.
	By the abstract duality, this is a tensor product of
	$U_{s,s+a}$ and $U_{s,s+b}$.
	
	For matroids $M,N$ of equal total rank on the same finite
	ground set $E$, the dual-rank formula
	\[
	r_{M^*}(X)=|X|-r_M(E)+r_M(E\setminus X)
	\qquad(X\subseteq E)
	\]
	shows that $M\preceq N$ implies $M^*\preceq N^*$.
	Since duality also preserves distinctness,
	\eqref{eq:duality} and Theorem~\ref{thm:birigidity}
	give~\eqref{eq:main}.
	
	By Lemma~\ref{lem:relaxation}, the relaxed circuit-hyperplane
	in the dual is the complement of $E(G_k^{a,b})$.
	Its cardinality is
	\[
	mn-r_{a,b}(m,n)=(m-a)(n-b)=s^2.
	\]
	Every edge incident with an added vertex belongs to
	$G_k^{a,b}$, so all edges of this complement lie between
	the original vertex classes of $G_k$.
\end{proof}

\begin{corollary}
	For all integers $a,b,k\ge2$, the matroids $N_{a,b,k}$
	and $T_{a,b,k}$ are non-representable over every field.
\end{corollary}

\begin{proof}
	By the author's theorem~\cite[Theorem~1.2]{Tyomkyn},
	every representable abstract $(a,b)$-birigidity matroid
	is less free than the characteristic-zero generic matroid
	with the same parameters. This contradicts
	Theorem~\ref{thm:birigidity} if $N_{a,b,k}$ is representable.
	Representability is preserved by duality, so the same
	holds for $T_{a,b,k}$.
\end{proof}

\begin{remark}
	One may relax all part-preserving permutation copies of
	$G_k^{a,b}$, obtaining a matroid invariant under permutations
	of either vertex class. Distinct original circuit-hyperplanes
	remain circuit-hyperplanes after any one of them is relaxed,
	so these relaxations can be performed successively.
	When $a=b$, the resulting matroid is also invariant
	under interchanging the two vertex classes.
\end{remark}

\begin{remark}
	The construction disproves the asserted maximality of the
	generic tensor product. It does not exclude a different
	freest tensor product of the same uniform matroids.
\end{remark}

\section{Hyperconnectivity and the Jackson--Tanigawa conjecture}\label{sec:hypercon}
The matroids $\HH_d(K_n)$ are symmetric under vertex permutations and compatible under restriction between different $n$. We therefore may, slightly abusing notation, speak of $\HH_d$-independent graphs and $\HH_d$-circuits.
On a bipartite graph, ordering the left vertices before the
right ones gives the signed $(d,d)$-birigidity representation. The next lemma relates their rank functions under edge contraction.

\begin{lemma}\label{lem:bipartite-contraction}
	Let $G\subseteq K_{m,n}$, where $m,n\ge3$, and let
	$ab\in E(G)$. Then
	\[
	r_{\HH_2}(G/ab)=r_{\BR_{2,2}(m,n)}(G)-1.
	\]
	Moreover, if $G$ is a circuit in $\BR_{2,2}(m,n)$,
	then $G/ab$ is a circuit in $\HH_2$.
\end{lemma}

\begin{proof}
	Take generic vertex vectors subject only to $p_a=p_b=p$.
	The resulting representation of $K_{m,n}$ still realizes
	the generic birigidity matroid: independent invertible
	changes of coordinates on the two vertex classes preserve
	all bipartite dependencies~\cite[Proposition~3.10(1)]{CR}
	and can be used to identify the two specified vectors.
	
	Replace the row blocks at $a$ and $b$ by
	a single block, their entry-wise sum, indexed by the merged vertex $c$. The column of $ab$ becomes zero, and the remaining
	columns give the generic hyperconnectivity representation
	of $K_{m,n}/ab$.
	
	Motions of the merged matrix correspond bijectively
	to motions of the original matrix satisfying
	$\dot p_a=\dot p_b$: the correspondence assigns the
	motion vector at $c$ to both $a$ and $b$ and leaves
	all other motion vectors unchanged.
	The original motions have the form $\dot p_u=Ap_u$
	on the left and $\dot p_v=A^\T p_v$ on the right,
	where $A\in\FF^{2\times2}$.
	Thus the equality becomes $Ap=A^\T p$.
	Since $p\ne0$, this holds precisely when $A$ is symmetric.
	Merging therefore removes two rows and reduces the
	motion-space dimension from $4$ to $3$, so the overall rank
	decreases by exactly $1$.
	
On the column span of $G$, the merging map has a
nonzero kernel, since the column of $ab$ is sent to $0$.
This kernel is one-dimensional, since the overall rank decreases by $1$. Hence the rank of $G$ also drops by exactly $1$.
	
If $G$ is a circuit, take its full-support stress.
Adding the row blocks at $a$ and $b$ preserves the
corresponding linear relation among the columns.
Deleting the zero column of $ab$ leaves a stress
of $G/ab$ with all remaining coefficients nonzero.
Since both the number of columns and the rank decrease
by $1$, the stress space remains one-dimensional.
Hence $G/ab$ is a circuit.
\end{proof}
\begin{lemma}\label{lem:contraction-closure}
	Let $D$ be a graph and let $c\in V(D)$ such that
	$D-c$ is bipartite, with vertex classes $L,R$.
	If distinct vertices $i,j$ belong to the same class
	and $ij\in\cl_{\HH_2}(D)$, then
	\[
	ci\in\cl_{\HH_2}(D)
	\qquad\text{and}\qquad
	cj\in\cl_{\HH_2}(D).
	\]
\end{lemma}

\begin{proof}
	For vertex vectors $p=(p_v)$ and motion vectors
	$\dot p=(\dot p_v)$, write
	\[
	\varepsilon_{p,\dot p}(uv)
	=p_v^\T\dot p_u-p_u^\T\dot p_v.
	\]
	Whenever $p$ represents $\HH_2$, an edge belongs to
	$\cl_{\HH_2}(D)$ precisely when this evaluation
	vanishes on every motion of $D$.
	
	By symmetry, suppose that $i,j\in L$.
	Choose generic vertex vectors $p$ and apply a common
	invertible change of coordinates to normalize
	\[
	p_c=(1,0)^\T
	\qquad\text{and}\qquad
	p_i=(0,1)^\T.
	\]
	The remaining vertex coordinates are algebraically
	independent.
	
	To prove that $ci\in\cl_{\HH_2}(D)$, it suffices to
	show that $\varepsilon_{p,\dot p}(ic)=0$ for every
	motion satisfying $\dot p_c=0$.
	Indeed, if a motion has value $(a,b)^\T$ at $c$,
	subtract the assignment $Sp_v$ at every vertex, where
	\[
	S=\begin{pmatrix}a&b\\b&0\end{pmatrix}.
	\]
	Since $S$ is symmetric, this changes no edge evaluation.
	At $c$, it subtracts $Sp_c=(a,b)^\T$, resulting in a $0$-vector.
	
	So, consider an arbitrary motion $\dot p$ of $D$ with
	$\dot p_c=0$. Put $J=\operatorname{diag}(1,-1)$
	and define a new $\HH_2$-representation and its motion as follows. Let
	\[
	p'_v=
	\begin{cases}
		Jp_v,&v\in R,\\
		p_v,&v\in L\cup\{c\},
	\end{cases}
	\qquad\text{and}\qquad
	\dot p'_v=
	\begin{cases}
		J\dot p_v,&v\in L,\\
		\dot p_v,&v\in R\cup\{c\}.
	\end{cases}
	\]
	The point coordinates have at most changed
	signs, so $p'$ also represents the generic matroid
	$\HH_2$.
	
	We check that $\dot p'$ is a motion of $D$ in the
	representation determined by $p'$.
	For an edge $uv$ with $u\in L$ and $v\in R$ we have
	\[
	\begin{aligned}
		\varepsilon_{p',\dot p'}(uv)
		&=(Jp_v)^\T J\dot p_u-p_u^\T\dot p_v=\varepsilon_{p,\dot p}(uv)=0.
	\end{aligned}
	\]
	For every edge $cv$ of $D$, its original equation
	reduces to $p_c^\T\dot p_v=0$, since $\dot p_c=0$.
	Thus the first coordinate of $\dot p_v$ is zero.
	The first coordinate of $\dot p'_v$ is the same,
	and $p'_c=p_c$ and $\dot p'_c=0$, so
	$\varepsilon_{p',\dot p'}(cv)=0$ as well.
	
    We can therefore apply the assumption
	$ij\in\cl_{\HH_2}(D)$ to both the motion $\dot p$
	in the representation $p$ and the motion $\dot p'$
	in the representation $p'$.
	
	Write
	\[
	p_j=(x,y)^\T,\qquad
	\dot p_i=(u,v)^\T
	\qquad\text{and}\qquad
	\dot p_j=(r,s)^\T.
	\]
	Since $i,j\in L$, their position vectors are unchanged:
	$p'_i=p_i$ and $p'_j=p_j$, whereas
	\[
	\dot p'_i=(u,-v)^\T
	\qquad\text{and}\qquad
	\dot p'_j=(r,-s)^\T.
	\]
	The two evaluations consequently give
	\[
	\varepsilon_{p,\dot p}(ij)=xu+yv-s=0
	\]
	and
	\[
	\varepsilon_{p',\dot p'}(ij)=xu-yv+s=0.
	\]
	Adding gives $2xu=0$. Since $x\ne0$ by genericity,
	we obtain $u=0$, and hence
	\[
	\varepsilon_{p,\dot p}(ic)
	=p_c^\T\dot p_i-p_i^\T\dot p_c=u=0.
	\]
	The motion $\dot p$ was arbitrary subject to
	$\dot p_c=0$, so the reduction above proves
	$ci\in\cl_{\HH_2}(D)$.
	Interchanging $i,j$ gives $cj\in\cl_{\HH_2}(D)$.
\end{proof}

For each integer $k\ge2$, put
\[
D_k=G_k/(a_0^1b_0^1).
\]
This creates no parallel edges, since the endpoints of an
edge in a bipartite graph have no common neighbour.
Thus $D_k$ has $8k+1$ vertices and $16k-1$ edges.

\begin{proposition}\label{prop:Dk}
	For every integer $k\ge2$, the graph $D_k$ is a
	circuit-hyperplane of $\HH_2(K_{8k+1})$.
\end{proposition}

\begin{proof}
	Write $a=a_0^1$, $b=b_0^1$, and let $c$ be the vertex
	obtained by identifying $a,b$. Put
	\[
	L'=L_k\setminus\{a\},\qquad
	R'=R_k\setminus\{b\}.
	\]
	By Proposition~\ref{prop:seed} and
	Lemma~\ref{lem:bipartite-contraction}, the graph $D_k$
	is a circuit in $\HH_2(K_{8k+1})$, with
	\[
	r_{\HH_2}(D_k)=16k-2.
	\]
	
Since $G_k$ is a flat of $\BR_{2,2}(4k+1,4k+1)$,
adding any missing edge from $L_k\times R_k$
increases its rank by $1$.
Applying Lemma~\ref{lem:bipartite-contraction}
before and after adding this edge shows that the
corresponding rank increase persists after contraction.
Thus every missing edge between $L'$ and $R'$,
and every missing edge incident with $c$, lies
outside the closure of $D_k$.

It remains to consider pairs of vertices in the same
class. By Lemma~\ref{lem:contraction-closure}, such
a pair can belong to the closure only if both
vertices are neighbours of $c$.
The neighbours of $c$ in $L'$ are
\[
\{a_1^1,x_0^1,x_0^2\}.
\]
We exclude all three pairs in this set using the fact
that every $K_4$ is a circuit: any five of its edges
span the remaining edge.

Suppose first that $x_0^1x_0^2$ belongs to the closure.
The four edges joining $x_0^1,x_0^2$ to $c,b_*$
already belong to $D_k$. Applying the $K_4$ circuit
property on
\[
\{c,b_*,x_0^1,x_0^2\}
\]
therefore puts $cb_*$ in the closure, contradicting
the fact that this missing edge lies outside it.

Now suppose that $a_1^1x_0^\sigma$ belongs to the closure
for some $\sigma\in\{1,2\}$. Interchanging $x_0^1,x_0^2$
preserves $D_k$ and, by the relabelling invariance of
generic hyperconnectivity, preserves its closure.
Thus both $a_1^1x_0^1$ and $a_1^1x_0^2$ belong to the
closure. The three edges joining $c$ to
$a_1^1,x_0^1,x_0^2$ already belong to $D_k$.
The $K_4$ circuit property on
\[
\{c,a_1^1,x_0^1,x_0^2\}
\]
now puts $x_0^1x_0^2$ in the closure, giving the
previous contradiction. 

For pairs within $R'$, the assertion follows symmetrically.
	
	Thus every missing edge lies outside the closure of
	$D_k$. Its rank is $16k-2$, whereas the ambient rank
	is $16k-1$, so it is a hyperplane. Together with the
	circuit assertion, this proves the proposition.
\end{proof}

Coning takes any circuit-hyperplane of $\HH_d(K_m)$
using every vertex to a circuit-hyperplane of
$\HH_{d+1}(K_{m+1})$.
To see this, let $d,m$ be integers with $d\ge1$ and
$m\ge d+1$, and let $G$ be a graph on $[m]$.
Add a new apex adjacent to every old vertex.
Normalize its vector to
$(0,\ldots,0,1)^\T\in\FF^{d+1}$, and write the old vectors
as $(q_i^\T,t_i)^\T$, where $q_i\in\FF^d$ and $i\in[m]$.
The coordinates of the $q_i$ and the scalars $t_i$
are generic. Put
\[
t=(t_1,\ldots,t_m)^{\T},
\]
and let $Q$ be the matrix with rows $q_i^{\T}$.

For an old stress $\omega$, define an $m\times m$ matrix
$\Omega$ by setting
\[
\Omega_{ij}=\omega(ij)
\qquad\text{and}\qquad
\Omega_{ji}=-\omega(ij)
\]
for every edge $ij$ with $i<j$, and setting all remaining
entries to $0$. Then $\Omega^\T=-\Omega$, and the stress
equation at vertex $i$ is
\[
\sum_j\Omega_{ij}q_j=0.
\]
Since $Q$ has rows $q_j^\T$, these equations are equivalent
to $\Omega Q=0$. Thus the stresses of $G$ correspond
exactly to the skew-symmetric matrices $\Omega$ supported
on $E(G)$ and satisfying $\Omega Q=0$.

The stresses of $\cone G$ in $\HH_{d+1}$ are represented
by precisely the following skew-symmetric matrices:
\begin{equation}\label{eq:hypercone}
	\widehat\Omega=
	\begin{pmatrix}
		\Omega&-\Omega t\\
		-t^{\T}\Omega&0
	\end{pmatrix},
	\qquad
	\Omega^{\T}=-\Omega,\quad \Omega Q=0.
\end{equation}
Indeed, the last-coordinate equations at the old vertices
force the new column to be $-\Omega t$; skew-symmetry
then determines the new row. The remaining equations
follow from $\Omega Q=0$ and $t^{\T}\Omega t=0$.
Thus restriction to the old edges is a linear bijection
between the stress spaces of $\cone G$ and $G$:
every old stress has the unique extension given
by~\eqref{eq:hypercone}. In particular, the two stress
spaces have the same dimension.
Since $m$ edges were added,
\[
r_{\HH_{d+1}}(\cone G)=r_{\HH_d}(G)+m.
\]

If $G$ is a circuit using every vertex, a nonzero old
stress has full support and no zero row. The independent scalars $t_i$ therefore make the stress
coefficient on every new edge nonzero, so $\cone G$ is also a circuit.
Every missing edge is old, and applying the rank formula
to $G$ and $G+f$ preserves the rank increase from adding
any missing edge $f$. Hence flatness is preserved.
Finally, for the overall rank of the matroids we have
\[
r_{d+1}(m+1)=r_d(m)+m.
\]
It follows that a circuit-hyperplane using every vertex
cones to a circuit-hyperplane in the next dimension.

\begin{proof}[Proof of Theorem~\ref{thm:hyper}]
	Put	$F_{d,k}=\cone^{d-2}(D_k)$.
	Proposition~\ref{prop:Dk} and the coning argument
	based on~\eqref{eq:hypercone} show that $F_{d,k}$ is
	a circuit-hyperplane with
	\[
	|E(F_{d,k})|=r_d(n)=\frac{d(16k+d-3)}2.
	\]
	Let $M_{d,k}$ be its relaxation in $\HH_d(K_n)$.
	For $d\ge3$, the graph $F_{d,k}$ contains triangles,
	so it is not a $K_{d+1,d+1}$; for $d=2$, its vertex
	count excludes equality with $K_{3,3}$.
	Its vertex count also excludes equality with $K_{d+2}$.
	Relaxation therefore preserves all circuits of both
	types and the total rank, and gives
	$\HH_d(K_n)\prec M_{d,k}$.
	
	Nguyen's characterization of abstract
	rigidity~\cite[Theorem~4.2]{CJT} says that rank $r_d(n)$
	together with all $K_{d+2}$ being circuits is sufficient.
	Finally, $n\ge2d+2$ is equivalent to $8k\ge d+3$.
\end{proof}

Bulavka and Tancer~\cite{BT} disproved the
Jackson--Tanigawa hyperconnectivity conjecture for
integers $d\ge3$ and $n\ge2d+2$.
Their explicit witness has rank $(d+1)^2$, strictly
below $r_d(n)$, and is incomparable with $\HH_d(K_n)$.
Their extension of this matroid to a maximal matroid for its prescribed circuits does not specify the rank of the latter.
Our construction gives a matroid strictly above
$\HH_d(K_n)$, with full rigidity rank, and includes
$d=2$.

\end{document}